\documentclass[a4paper,11pt]{amsart}
\usepackage{amsmath,amssymb,amsfonts,latexsym}

\newtheorem{theorem}{Theorem}[section]

\input{tcilatex}

\begin{document}
\title[$q$\textbf{-Dedekind-type D-C sums}]{\textbf{A note on the} $q$%
\textbf{-Dedekind-type Daehee-Changhee sums with weight }$\alpha $\textbf{\
arising from modified }$q$\textbf{-Genocchi polynomials with weight }$\alpha 
$}
\author[\textbf{S. Araci}]{\textbf{Serkan Araci}}
\address{\textbf{University of Gaziantep, Faculty of Science and Arts,
Department of Mathematics, 27310 Gaziantep, TURKEY}}
\email{\textbf{mtsrkn@hotmail.com; saraci88@yahoo.com.tr; mtsrkn@gmail.com}}
\author[\textbf{M. Acikgoz}]{\textbf{Mehmet Acikgoz}}
\address{\textbf{University of Gaziantep, Faculty of Science and Arts,
Department of Mathematics, 27310 Gaziantep, TURKEY}}
\email{\textbf{acikgoz@gantep.edu.tr}}
\author[\textbf{A. Esi}]{\textbf{Ayhan Esi}}
\address{\textbf{University of Adiyaman, Faculty of Science and Arts,
Department of Mathematics, 02040 Adiyaman, TURKEY}}
\email{\textbf{aesi23@hotmail.com}}

\begin{abstract}
In the present paper, our objective is to treat a $p$-adic continuous
function for an odd prime to inside a $p$-adic $q$-analogue of the higher
order Dedekind-type sums with weight $\alpha $ in connection with modified $%
q $-Genocchi polynomials with weight $\alpha $ by using $p$-adic invariant $%
q $-integral on $%
\mathbb{Z}
_{p}$.

\vspace{2mm}\noindent \textsc{2010 Mathematics Subject Classification.}
11S80, 11B68.

\vspace{2mm}

\noindent \textsc{Keywords and phrases.} Dedekind sums, $q$-Dedekind-type
sums, $p$-adic $q$-integral on $%
\mathbb{Z}
_{p}$, modified $q$-Genocchi polynomials with weight $\alpha .$
\end{abstract}

\thanks{}
\maketitle




\section{\textbf{Introduction}}


Assume that $p$ be a fixed odd prime number. We now begin with the
definition of the following notations. Let $%
\mathbb{Q}
_{p}$ be the field $p$-adic rational numbers and let $%
\mathbb{C}
_{p}$ be the completion of algebraic closure of $%
\mathbb{Q}
_{p}$.

Thus, 
\begin{equation*}
\boldsymbol{%
\mathbb{Q}
}_{p}=\left\{ x=\sum_{n=-k}^{\infty }a_{n}p^{n}:0\leq a_{n}<p\right\} \text{.%
}
\end{equation*}

Then $%
\mathbb{Z}
_{p}$ is integral domain, which is given by 
\begin{equation*}
\boldsymbol{%
\mathbb{Z}
}_{p}=\left\{ x=\sum_{n=0}^{\infty }a_{n}p^{n}:0\leq a_{n}<p\right\}
\end{equation*}%
or%
\begin{equation*}
\boldsymbol{%
\mathbb{Z}
}_{p}=\left\{ x\in 
\mathbb{Q}
_{p}:\left\vert x\right\vert _{p}\leq 1\right\} \text{.}
\end{equation*}

We suppose that $q\in 
\mathbb{C}
_{p}$ with $\left\vert 1-q\right\vert _{p}<1$ as an indeterminate. The $p$%
-adic absolute value $\left\vert .\right\vert _{p}$, is normally introduced
by 
\begin{equation*}
\left\vert x\right\vert _{p}=\frac{1}{p^{n}}
\end{equation*}%
where $x=p^{n}\frac{s}{t}$ with $\left( p,s\right) =\left( p,t\right)
=\left( s,t\right) =1$ and $n\in 
\mathbb{Q}
$ (see [1-15]).

The $p$-adic $q$-Haar distribution was originally introduced by Kim as
follows: For each postive integer $n$,%
\begin{equation*}
\mu _{q}\left( a+p^{n}%
\mathbb{Z}
_{p}\right) =\left( -q\right) ^{a}\frac{\left( 1+q\right) }{1+q^{p^{n}}}
\end{equation*}%
for $0\leq a<p^{n}$ and this can be extended to a measure on $%
\mathbb{Z}
_{p}$ (for details, see [1-7]).

In \cite{Araci}, modified $q$-Genocchi polynomials with weight $\left(
\alpha ,\beta \right) $ are defined by Araci \textit{et al.} as follows:%
\begin{equation}
\widetilde{G}_{n,q}^{\left( \alpha ,\beta \right) }\left( x\right) =n\int_{%
\mathbb{Z}
_{p}}q^{-\beta \xi }\left( \frac{1-q^{\alpha \left( x+\xi \right) }}{%
1-q^{\alpha }}\right) ^{n-1}d\mu _{q^{\beta }}\left( \xi \right)
\label{equation 1}
\end{equation}%
for $n\in 
\mathbb{Z}
_{+}:=\left\{ 0,1,2,3,\cdots \right\} $. We easily see that 
\begin{equation*}
\lim_{q\rightarrow 1}\widetilde{G}_{n,q}^{\left( \alpha \right) }\left(
x\right) =G_{n}\left( x\right)
\end{equation*}%
where $G_{n}\left( x\right) $ are Genocchi polynomials, which are given in
the form:%
\begin{equation*}
\sum_{n=0}^{\infty }G_{n}\left( x\right) \frac{t^{n}}{n!}=e^{tx}\frac{2t}{%
e^{t}+1},\text{ }\left\vert t\right\vert <\pi
\end{equation*}%
(for details, see \cite{Araci 2}). Taking $x=0$ into (\ref{equation 1}),
then we have $\widetilde{G}_{n,q}^{\left( \alpha ,\beta \right) }\left(
0\right) :=\widetilde{G}_{n,q}^{\left( \alpha ,\beta \right) }$ are called
modified $q$-Genocchi numbers with weight $\left( \alpha ,\beta \right) $.

It seems to be interesting for studying equation (\ref{equation 1}) at $%
\beta =1$. Then, we can state the following:%
\begin{equation}
\widetilde{G}_{n,q}^{\left( \alpha ,1\right) }\left( x\right) :=\widetilde{G}%
_{n,q}^{\left( \alpha \right) }\left( x\right) =n\int_{%
\mathbb{Z}
_{p}}q^{-\xi }\left( \frac{1-q^{\alpha \left( x+\xi \right) }}{1-q^{\alpha }}%
\right) ^{n-1}d\mu _{q}\left( \xi \right) \text{.}  \label{equation 9}
\end{equation}%
where $\widetilde{G}_{n,q}^{\left( \alpha \right) }\left( x\right) $ are
called modified $q$-Genocchi polynomials with weight $\alpha $.

Modified $q$-Genocchi numbers and polynomials with weight $\alpha $ have the
following identities:%
\begin{eqnarray}
\widetilde{G}_{n+1,q}^{\left( \alpha \right) } &=&\left( n+1\right) \frac{1+q%
}{\left( 1-q^{\alpha }\right) ^{n}}\sum_{l=0}^{n}\binom{n}{l}\left(
-1\right) ^{l}\frac{1}{1+q^{\alpha l}}\text{,}  \label{equation 2} \\
\widetilde{G}_{n+1,q}^{\left( \alpha \right) }\left( x\right) &=&\left(
n+1\right) \frac{1+q}{\left( 1-q^{\alpha }\right) ^{n}}\sum_{l=0}^{n}\binom{n%
}{l}\left( -1\right) ^{l}\frac{q^{\alpha lx}}{1+q^{\alpha l}}\text{,}
\label{equation 3} \\
\widetilde{G}_{n,q}^{\left( \alpha \right) }\left( x\right) &=&q^{-\alpha
x}\sum_{l=0}^{n}\binom{n}{l}q^{\alpha lx}\widetilde{G}_{l,q}^{\left( \alpha
\right) }\left( \frac{1-q^{\alpha x}}{1-q^{\alpha }}\right) ^{n-l}\text{.}
\label{equation 4}
\end{eqnarray}

Additionally, for $d$ odd natural number, we have 
\begin{equation}
\widetilde{G}_{n,q}^{\left( \alpha \right) }\left( dx\right) =\left( \frac{%
1+q}{1+q^{d}}\right) \left( \frac{1-q^{\alpha d}}{1-q^{\alpha }}\right)
^{n-1}\sum_{a=0}^{d-1}\left( -1\right) ^{a}\widetilde{G}_{n,q}^{\left(
\alpha \right) }\left( x+\frac{a}{d}\right) \text{,}  \label{equation 5}
\end{equation}%
(for details about this subject, see \cite{Araci}).

For any positive integer $h,k$ and $m$, Dedekind-type D-C sums are given by
Kim in \cite{Kim 1}, \cite{Kim 2} and \cite{Kim 3} as follows:%
\begin{equation*}
S_{m}\left( h,k\right) =\sum_{M=1}^{k-1}\left( -1\right) ^{M-1}\frac{M}{k}%
\overline{E}_{m}\left( \frac{hM}{k}\right)
\end{equation*}%
where $\overline{E}_{m}\left( x\right) $ are the $m$-th periodic Euler
function.

In 2011, Taekyun Kim introduced weighted $q$-Bernoulli numbers and
polynomials in \cite{Kim 8}. He derived not only new but also ineteresting
properties for weighted $q$-Bernoulli numbers and polynomials. In \cite%
{Araci}, Araci \textit{et al}. extended Kim's method for $q$-Genocchi
polynomials and also they defined $q$-Genocchi numbers and polynomials with
weight ($\alpha ,\beta $).

In \cite{Kim 2}, Kim has given some fascinating properties for Dedekind-type
D-C sums. He firstly considered a $p$-adic continuous function for an odd
prime number to contain a $p$-adic $q$-analogue of the higher order
Dedekind-type D-C sums $k^{m}S_{m+1}\left( h,k\right) $. In previous paper 
\cite{Araci 5}, Araci and Acikgoz also introduced the definition of the
extended $q$-Dedekind-type sums and given the relation between extended $q$%
-Euler polynomials.

By using $p$-adic invariant $q$-integral on $%
\mathbb{Z}
_{p}$, in this paper, we shall give the definition of $q$-Dedekind-type sums
with weight $\alpha $. Also, we shall derive interesting property for $q$%
-Dedekind sums with weight $\alpha $ arising from modified $q$-Genocchi
polynomials with weight $\alpha $.

\section{$q$\textbf{-Dedekind-type D-C sums with weight }$\protect\alpha $%
\textbf{\ related to modified }$q$\textbf{-Genocchi polynomials with weight }%
$\protect\alpha $}

If $x$ is a $p$-adic integer, then $w\left( x\right) $ is the unique
solution of $w\left( x\right) =w\left( x\right) ^{p}$ that is congruent to $x%
\func{mod}p$. It can also be defined by%
\begin{equation*}
w\left( x\right) =\lim_{n\rightarrow \infty }x^{p^{n}}\text{.}
\end{equation*}

The multiplicative group of $p$-adic units is a product of the finite group
of roots of unity, and a group isomorphic to the $p$-adic integers. The
finite group is cylic of order $p-1$ or $2$, as $p$ is odd or even,
respectively, and so it is isomorphic. Actually, the teichm\"{u}ller
character gives a canonical isomorphism between these two groups.

Let $w$ be the $Teichm\ddot{u}ller$ character ($\func{mod}p$). For $x\in 
\mathbb{Z}
_{p}^{\ast }$ $:=%
\mathbb{Z}
_{p}/p%
\mathbb{Z}
_{p}$, set%
\begin{equation*}
\left\langle x:q\right\rangle =w^{-1}\left( x\right) \left( \frac{1-q^{x}}{%
1-q}\right) \text{.}
\end{equation*}

Let $a$ and $N$ be positive integers with $\left( p,a\right) =1$ and $p\mid
N $. We now introduce the following%
\begin{equation*}
\widetilde{A}_{q}^{\left( \alpha \right) }\left( s,a,N:q^{N}\right)
=w^{-1}\left( a\right) \left\langle a:q^{\alpha }\right\rangle
^{s}q^{-\alpha a}\sum_{j=0}^{\infty }\binom{s}{j}q^{\alpha aj}\left( \frac{%
1-q^{\alpha N}}{1-q^{\alpha a}}\right) ^{j}\widetilde{G}_{j,q^{N}}^{\left(
\alpha \right) }\text{.}
\end{equation*}

Obviously, if $m+1\equiv 0(\func{mod}p-1)$, then%
\begin{eqnarray*}
\widetilde{A}_{q}^{\left( \alpha \right) }\left( m,a,N:q^{N}\right)
&=&\left( \frac{1-q^{\alpha a}}{1-q^{\alpha }}\right) ^{m}q^{-\alpha
a}\sum_{j=0}^{m}\binom{m}{j}q^{\alpha aj}\widetilde{G}_{j,q^{N}}^{\left(
\alpha \right) }\left( \frac{1-q^{\alpha N}}{1-q^{\alpha a}}\right) ^{j} \\
&=&\left( \frac{1-q^{\alpha N}}{1-q^{\alpha }}\right) ^{m}\int_{%
\mathbb{Z}
_{p}}q^{-N\xi }\left( \frac{1-q^{\alpha N\left( \xi +\frac{a}{N}\right) }}{%
1-q^{\alpha N}}\right) ^{m}d\mu _{q^{N}}\left( \xi \right) \text{.}
\end{eqnarray*}

Then, $\widetilde{A}_{q}^{\left( \alpha \right) }\left( m,a,N:q^{N}\right) $
is a continuous $p$-adic analogue of 
\begin{equation*}
\left( \frac{1-q^{\alpha N}}{1-q^{\alpha }}\right) ^{m}\frac{\widetilde{G}%
_{m+1,q^{N}}^{\left( \alpha \right) }\left( \frac{a}{N}\right) }{m+1}\text{.}
\end{equation*}

Let $\left[ .\right] $ be the Gauss' symbol and let $\left\{ x\right\} =x-%
\left[ x\right] $. Thus, we are now ready to treat $q$-analogue of the
higher order Dedekind-type D-C sums $\widetilde{S}_{m,q}^{\left( \alpha
\right) }\left( h,k:q^{l}\right) $ in the form: 
\begin{equation*}
\widetilde{Y}_{m,q}^{\left( \alpha \right) }\left( h,k:q^{l}\right)
=\sum_{M=1}^{k-1}\left( -1\right) ^{M-1}\left( \frac{1-q^{\alpha M}}{%
1-q^{\alpha k}}\right) \int_{%
\mathbb{Z}
_{p}}q^{-l\xi }\left( \frac{1-q^{\alpha \left( l\xi +l\left\{ \frac{hM}{k}%
\right\} \right) }}{1-q^{\alpha l}}\right) ^{m}d\mu _{q^{l}}\left( \xi
\right) \text{.}
\end{equation*}

If $m+1\equiv 0\left( \func{mod}p-1\right) $%
\begin{eqnarray*}
&&\left( \frac{1-q^{\alpha k}}{1-q^{\alpha }}\right)
^{m+1}\sum_{M=1}^{k-1}\left( -1\right) ^{M-1}\left( \frac{1-q^{\alpha M}}{%
1-q^{\alpha k}}\right) \int_{%
\mathbb{Z}
_{p}}q^{-\xi k}\left( \frac{1-q^{\alpha k\left( \xi +\frac{hM}{k}\right) }}{%
1-q^{\alpha k}}\right) ^{m}d\mu _{q^{k}}\left( \xi \right) \\
&=&\sum_{M=1}^{k-1}\left( -1\right) ^{M-1}\left( \frac{1-q^{\alpha M}}{%
1-q^{\alpha }}\right) \left( \frac{1-q^{\alpha k}}{1-q^{\alpha }}\right)
^{m}\int_{%
\mathbb{Z}
_{p}}q^{-\xi k}\left( \frac{1-q^{\alpha k\left( \xi +\frac{hM}{k}\right) }}{%
1-q^{\alpha k}}\right) ^{m}d\mu _{q^{k}}\left( \xi \right)
\end{eqnarray*}%
where $p\mid k$, $\left( hM,p\right) =1$ for each $M$. By means of the
equation (\ref{equation 1}), we easily derive the following: 
\begin{gather}
\left( \frac{1-q^{\alpha k}}{1-q^{\alpha }}\right) ^{m+1}\widetilde{Y}%
_{m,q}^{\left( \alpha \right) }\left( h,k:q^{k}\right)  \label{equation 6} \\
=\sum_{M=1}^{k-1}\left( \frac{1-q^{\alpha M}}{1-q^{\alpha }}\right) \left( 
\frac{1-q^{\alpha k}}{1-q^{\alpha }}\right) ^{m}\left( -1\right) ^{M-1}\int_{%
\mathbb{Z}
_{p}}q^{-\xi k}\left( \frac{1-q^{\alpha k\left( \xi +\frac{hM}{k}\right) }}{%
1-q^{\alpha k}}\right) ^{m}d\mu _{q^{k}}\left( \xi \right)  \notag \\
=\sum_{M=1}^{k-1}\left( -1\right) ^{M-1}\left( \frac{1-q^{\alpha M}}{%
1-q^{\alpha }}\right) \widetilde{A}_{q}^{\left( \alpha \right) }\left(
m,\left( hM\right) _{k}:q^{k}\right)  \notag
\end{gather}%
where $(hM)_{k}$ denotes the integer $x$ such that $0\leq x<n$ and $x\equiv
\alpha \left( \func{mod}k\right) $.

It is easy to indicate the following: 
\begin{gather}
\int_{%
\mathbb{Z}
_{p}}q^{-\xi }\left( \frac{1-q^{\alpha \left( x+\xi \right) }}{1-q^{\alpha }}%
\right) ^{k}d\mu _{q}\left( \xi \right)  \label{equation 7} \\
=\left( \frac{1-q^{\alpha m}}{1-q^{\alpha }}\right) ^{k}\frac{1+q}{1+q^{m}}%
\sum_{i=0}^{m-1}\left( -1\right) ^{i}\int_{%
\mathbb{Z}
_{p}}q^{-m\xi }\left( \frac{1-q^{\alpha m\left( \xi +\frac{x+i}{m}\right) }}{%
1-q^{\alpha m}}\right) ^{k}d\mu _{q^{m}}\left( \xi \right) \text{.}  \notag
\end{gather}

Due to (\ref{equation 6}) and (\ref{equation 7}), we get%
\begin{gather}
\left( \frac{1-q^{\alpha N}}{1-q^{\alpha }}\right) ^{m}\int_{%
\mathbb{Z}
_{p}}q^{-\xi N}\left( \frac{1-q^{\alpha N\left( \xi +\frac{a}{N}\right) }}{%
1-q^{\alpha N}}\right) ^{m}d\mu _{q^{N}}\left( \xi \right)
\label{equation 8} \\
=\frac{1+q^{N}}{1+q^{Np}}\sum_{i=0}^{p-1}\left( -1\right) ^{i}\left( \frac{%
1-q^{\alpha Np}}{1-q^{\alpha }}\right) ^{m}\int_{%
\mathbb{Z}
_{p}}q^{-\xi pN}\left( \frac{1-q^{\alpha pN\left( \xi +\frac{a+iN}{pN}%
\right) }}{1-q^{\alpha pN}}\right) ^{m}d\mu _{q^{pN}}\left( \xi \right) 
\text{.}  \notag
\end{gather}

Via the (\ref{equation 6}), (\ref{equation 7}) and (\ref{equation 8}), we
discover the following $p$-adic integration: 
\begin{equation*}
\widetilde{A}_{q}^{\left( \alpha \right) }\left( s,a,N:q^{N}\right) =\frac{%
1+q^{N}}{1+q^{Np}}\sum_{\underset{a+iN\neq 0(\func{mod}p)}{0\leq i\leq p-1}%
}\left( -1\right) ^{i}\widetilde{A}_{q}^{\left( \alpha \right) }\left(
s,\left( a+iN\right) _{pN},p^{N}:q^{pN}\right) \text{.}
\end{equation*}

In the other words,%
\begin{gather*}
\widetilde{A}_{q}^{\left( \alpha \right) }\left( m,a,N:q^{N}\right) =\left( 
\frac{1-q^{\alpha N}}{1-q^{\alpha }}\right) ^{m}\int_{%
\mathbb{Z}
_{p}}q^{-\xi N}\left( \frac{1-q^{\alpha N\left( \xi +\frac{a}{N}\right) }}{%
1-q^{\alpha N}}\right) ^{m}d\mu _{q^{N}}\left( \xi \right) \\
-\left( \frac{1-q^{\alpha Np}}{1-q^{\alpha }}\right) ^{m}\int_{%
\mathbb{Z}
_{p}}q^{-\xi N}\left( \frac{1-q^{\alpha pN\left( \xi +\frac{a+iN}{pN}\right)
}}{1-q^{\alpha pN}}\right) ^{m}d\mu _{q^{pN}}\left( \xi \right)
\end{gather*}%
where $\left( p^{-1}a\right) _{N}$ denotes the integer $x$ with $0\leq x<N$, 
$px\equiv a\left( \func{mod}N\right) $ and $m$ is integer with $m+1\equiv 0(%
\func{mod}p-1)$. Then we derive the following%
\begin{gather*}
\sum_{M=1}^{k-1}\left( -1\right) ^{M-1}\left( \frac{1-q^{\alpha M}}{%
1-q^{\alpha }}\right) \widetilde{A}_{q}^{\left( \alpha \right) }\left(
m,hM,k:q^{k}\right) \\
=\left( \frac{1-q^{\alpha k}}{1-q^{\alpha }}\right) ^{m+1}\widetilde{Y}%
_{m,q}^{\left( \alpha \right) }\left( h,k:q^{k}\right) -\left( \frac{%
1-q^{\alpha k}}{1-q^{\alpha }}\right) ^{m+1}\left( \frac{1-q^{\alpha kp}}{%
1-q^{\alpha k}}\right) \widetilde{Y}_{m,q}^{\left( \alpha \right) }\left(
\left( p^{-1}h\right) ,k:q^{pk}\right)
\end{gather*}%
where $p\nmid k$ and $p\nmid hm$ for each $M$. Thus, we get the following
definition.

\begin{definition}
Let $h,k$ be positive integer with $\left( h,k\right) =1$, $p\nmid k$. For $%
s\in 
\mathbb{Z}
_{p},$ we define $p$-adic Dedekind-type DC sums as follows:  
\begin{equation*}
\widetilde{Y}_{p,q}^{\left( \alpha \right) }\left( s:h,k:q^{k}\right)
=\sum_{M=1}^{k-1}\left( -1\right) ^{M-1}\left( \frac{1-q^{\alpha M}}{%
1-q^{\alpha }}\right) \widetilde{A}_{q}^{\left( \alpha \right) }\left(
m,hM,k:q^{k}\right) \text{.}
\end{equation*}
\end{definition}

As a result of the above definition, we derive the following theorem.

\begin{theorem}
For $m+1\equiv 0(\func{mod}p-1)$ and $\left( p^{-1}a\right) _{N}$ denotes
the integer $x$ with $0\leq x<N$, $px\equiv a\left( \func{mod}N\right) $,
then we have 
\begin{gather*}
\widetilde{Y}_{p,q}^{\left( \alpha \right) }\left( s:h,k:q^{k}\right)
=\left( \frac{1-q^{\alpha k}}{1-q^{\alpha }}\right) ^{m+1}\widetilde{Y}%
_{m,q}^{\left( \alpha \right) }\left( h,k:q^{k}\right) \\
-\left( \frac{1-q^{\alpha k}}{1-q^{\alpha }}\right) ^{m+1}\left( \frac{%
1-q^{\alpha kp}}{1-q^{\alpha k}}\right) \widetilde{Y}_{m,q}^{\left( \alpha
\right) }\left( \left( p^{-1}h\right) ,k:q^{pk}\right) \text{.}
\end{gather*}
\end{theorem}




\begin{thebibliography}{99}
\bibitem{Kim 1} T. Kim, A note on $p$-adic $q$-Dedekind sums, \textit{C. R.
Acad. Bulgare Sci.} 54 (\textbf{2001}), 37--42.

\bibitem{Kim 2} T. Kim, Note on $q$-Dedekind-type sums related to $q$-Euler
polynomials, \textit{Glasgow Math. J.} 54 (\textbf{2012}), 121-125.

\bibitem{Kim 3} T. Kim, Note on Dedekind type DC sums, \textit{Adv. Stud.
Contemp. Math}. 18 (\textbf{2009}), 249--260.

\bibitem{Kim 4} T. Kim, The modified $q$-Euler numbers and polynomials, 
\textit{Adv. Stud. Contemp. Math.} 16 (\textbf{2008}), 161--170.

\bibitem{Kim 5} T. Kim, $q$-Volkenborn integration, \textit{Russ. J. Math.
Phys.} 9 (\textbf{2002}), 288--299.

\bibitem{Kim 6} T. Kim, On $p$-adic interpolating function for $q$-Euler
numbers and its derivatives, \textit{J. Math. Anal. Appl.} 339 (\textbf{2008}%
), 598--608.

\bibitem{Kim 7} T. Kim, On a $q$-analogue of the $p$-adic log gamma
functions and related integrals, \textit{J. Number Theory} 76 (\textbf{1999}%
), 320-329.

\bibitem{Kim 8} T. Kim, On the weighted $q$-Bernoulli numbers and
polynomials, Advanced Studies in Contemporary Mathematics, vol. 21, no. 2,
pp. 207--215, 2011.

\bibitem{Kim 9} T. Kim and J. Choi, On the $q$-Euler Numbers and Polynomials
with Weight $0$, \textit{Abstract and Applied Analysis}, vol. \textbf{2012},
Article ID 795304, 7 pages.

\bibitem{Simsek} Y. Simsek, $q$-Dedekind type sums related to $q$-zeta
function and basic $L$-series, \textit{J. Math. Anal. Appl.} 318 (\textbf{%
2006}), 333-351.

\bibitem{Acikgoz} M. Acikgoz, Y. Simsek, On multiple interpolation function
of the N\"{o}rlund-type $q$-Euler polynomials, \textit{Abst. Appl. Anal. }%
2009 (\textbf{2009}), Article ID 382574, 14 pages.

\bibitem{Araci} S. Araci, M. Acikgoz, Feng Qi and H. Jolany, A note on the
modified $q$-Genocchi numbers and polynomials with weight $\left( \alpha
,\beta \right) $ and their interpolation function at negative integers, 
\textit{Fasc. Math. Journal (}accepted for publication\textit{).}

\bibitem{Araci 1} S. Araci, M. Acikgoz and K. H. Park, A note on the $q$%
-analogue of Kim's $p$-adic $\log $ gamma type functions associated with $q$%
-extension of Genocchi and Euler numbers with weight $\alpha $, accepted in 
\textit{Bulletin of the Korean Mathematical Society}.

\bibitem{Araci 2} S. Araci, M. Acikgoz and J. J. Seo, A study on the
weighted $q$-Genocchi numbers and polynomials with their interpolation
function, \textit{Honam Mathematical J}. 34 (\textbf{2012}), No. 1, pp.
11-18.

\bibitem{Araci 3} S. Araci, D. Erdal and J. J. Seo, A study on the fermionic 
$p$-adic $q$-integral representation on $%
\mathbb{Z}
_{p}$ associated with weighted $q$-Bernstein and $q$-Genocchi polynomials, 
\textit{Abstract and Applied Analysis}, Volume \textbf{2011}, Article ID
649248, 10 pages.

\bibitem{Araci 4} S. Araci, M. Acikgoz and J. J. Seo, Explicit formulas
involving $q$-Euler numbers and polynomials, \textit{Abstract and Applied
Analysis}, Volume \textbf{2012}, Article ID 298531, 11 pages.

\bibitem{Araci 5} S. Araci, and M. Acikgoz, Extended $q$-Dedekind-type sums
Daehee-Changhee sums associated with extended $q$-Euler polynomials, \textit{%
submitted.}
\end{thebibliography}
\end{document}